\documentclass[11pt]{amsart}

\usepackage[margin=1.15in]{geometry}
\usepackage{amsmath,amssymb,amsthm,mathtools}
\usepackage{enumitem}
\usepackage{booktabs}
\usepackage{array}
\usepackage[hidelinks]{hyperref}

\numberwithin{equation}{section}
\theoremstyle{plain}
\newtheorem{theorem}{Theorem}[section]
\newtheorem{proposition}[theorem]{Proposition}
\newtheorem{lemma}[theorem]{Lemma}
\newtheorem{corollary}[theorem]{Corollary}
\newtheorem{thmx}{Theorem}

\newtheorem{corx}[thmx]{Corollary}

\theoremstyle{definition}
\newtheorem{definition}[theorem]{Definition}

\theoremstyle{remark}
\newtheorem{remark}[theorem]{Remark}

\newcommand{\Z}{\mathbb{Z}}
\newcommand{\C}{\mathbb{C}}
\newcommand{\PP}{\mathbb{P}}
\newcommand{\cO}{\mathcal{O}}

\newcommand{\cB}{\mathcal{B}}
\newcommand{\cC}{\mathcal{C}}
\newcommand{\cD}{\mathcal{D}}
\newcommand{\Db}{\mathrm{D}^{\mathrm{b}}}
\newcommand{\Pset}{\mathsf{P}}
\newcommand{\aord}{\alpha_{\mathrm{ord}}}
\newcommand{\acomp}{\alpha_{\mathrm{C}}}
\newcommand{\gen}[1]{\langle #1\rangle}
\newcommand{\fto}{\longrightarrow}
\newcommand{\rev}{\mathrm{rev}}
\newcommand{\sm}{\setminus}
\DeclareMathOperator{\Fl}{Fl}
\DeclareMathOperator{\Sym}{Sym}
\DeclareMathOperator{\Pic}{Pic}
\DeclareMathOperator{\rk}{rk}
\DeclareMathOperator{\Wing}{Wing}
\DeclareMathOperator{\Ext}{Ext}
\DeclareMathOperator{\Hom}{Hom}
\DeclareMathOperator{\id}{id}

\begin{document}

\title{Exceptional sequences of line bundles on $\Fl(1,\ell;\ell+1)$}

\author{Likai Yang}
\address{Department of Mathematics, Michigan State University, East Lansing, MI 48824, USA}
\email{yanglika@msu.edu}


\begin{abstract}
Let $X_\ell=\Fl(1,\ell;\ell+1)\subset\PP^\ell\times(\PP^\ell)^\vee$ be the incidence variety of points and hyperplanes in $\PP^\ell$. We prove that for every $\ell\ge 3$, every exceptional sequence of line bundles on $X_\ell$ of maximal length $\ell(\ell+1)=\rk K_0(X_\ell)$ is full. Together with the case $\ell=2$, due to Altmann, Hochenegger and Witt, this shows that for all $\ell\ge2$ an exceptional sequence of line bundles on $X_\ell$ is full if and only if it has maximal length. The proof translates exceptionality into a combinatorial condition on lattice points of $\Pic(X_\ell)\cong\Z^2$ and reduces fullness to sequences containing a non-orthogonal descent. A capacity count for admissibly ordered subsets rules out all but one family of descents; for that family we classify the possible configurations and verify generation using Euler sequences and Koszul complexes.
\end{abstract}

\maketitle

\section{Introduction}\label{sec:intro}

Let $X$ be a smooth projective variety over $\C$ and $\Db(X)$ be its bounded derived category of coherent sheaves. An object $E\in\Db(X)$ is \emph{exceptional} if $\Hom(E,E)=\C$ and $\Ext^i(E,E)=0$ for $i\neq0$. A sequence $(E_1,\dots,E_n)$ of exceptional objects is an \emph{exceptional sequence} if $\Ext^\bullet(E_j,E_i)=0$ for all $j>i$, and it is \emph{full} if $E_1,\dots,E_n$ generate $\Db(X)$ as a triangulated category. The Euler pairing is unitriangular on the classes $[E_1],\dots,[E_n]$, so these classes are linearly independent. Hence, if $K_0(X)$ is free of finite rank (for instance, if $\Db(X)$ has a full exceptional sequence), every exceptional sequence has length at most $\rk K_0(X)$, and full exceptional sequences attain this bound. We say that an exceptional sequence has \emph{maximal length} if its length equals $\rk K_0(X)$.

The converse question, whether an exceptional sequence of maximal length must be full, is subtle. If such a sequence is not full, its right orthogonal is a nonzero admissible subcategory with torsion Grothendieck group. If, moreover, $X$ admits some full exceptional sequence, this subcategory has vanishing Grothendieck group and vanishing Hochschild homology, i.e.\ it is a phantom category. Phantom categories do exist~\cite{GO13,BGKS15}. Kuznetsov conjectured that phantoms cannot occur inside a category that is generated by an exceptional sequence, i.e.\ that in such a category every exceptional sequence of maximal length is full~\cite[Conjecture~1.10]{Kuz14}. Krah disproved this conjecture: the blow-up of $\PP^2$ in ten general points carries an exceptional sequence of line bundles of maximal length that is not full~\cite{Kra24}. On the other hand, fullness of maximal exceptional sequences of \emph{line bundles} is known in several important cases: for del Pezzo surfaces~\cite{EL16} (building on the toric systems of Hille and Perling~\cite{HP11}), for smooth projective toric varieties of Picard rank two~\cite{AW24}, and for the flag variety $\Fl(1,2;3)$~\cite{AHW25}. We also recall that on toric varieties exceptional sequences of line bundles may fail to reach the maximal length altogether~\cite{Efi14}.

In this paper we consider the point--hyperplane incidence variety
\[
X_\ell=\Fl(1,\ell;\ell+1)=\{(\lambda,H)\in\PP^\ell\times(\PP^\ell)^\vee:\ \lambda\in H\},\qquad \ell\ge2 .
\]
It is a rational homogeneous variety $\mathrm{SL}_{\ell+1}/P$ of dimension $2\ell-1$ and Picard rank two, and each of the two projections makes it a $\PP^{\ell-1}$-bundle over $\PP^\ell$; in particular $\rk K_0(X_\ell)=\ell(\ell+1)$. For $\ell\ge2$ it is not toric. Altmann, Hochenegger and Witt studied exceptional sequences of line bundles on projectivizations of vector bundles on $\PP^\ell$, with special attention to $X_\ell=\PP(\Omega_{\PP^\ell}(1))$; for $\ell=2$ they classified these sequences and proved that every one of maximal length is full~\cite[Theorem~C]{AHW25}. For $\ell\ge3$ the question remained open. Our main result settles it.

\begin{thmx}\label{thm:main}
Let $\ell\ge3$. Every exceptional sequence of $\ell(\ell+1)$ line bundles on $X_\ell$ is full.
\end{thmx}

Combining Theorem~\ref{thm:main} with~\cite[Theorem~C]{AHW25}, we obtain the follow corollary.

\begin{corx}\label{cor:main}
Let $\ell\ge2$ and let $L_1,\dots,L_n$ be line bundles on $X_\ell$. Then $(L_1,\dots,L_n)$ is a full exceptional sequence if and only if it is an exceptional sequence and $n=\ell(\ell+1)$.
\end{corx}

\subsection*{Lattice model.}
We identify the line bundle $L(a,b)=\cO(a)\boxtimes\cO(b)|_{X_\ell}$ with the lattice point $(a,b)\in\Z^2$. By the Borel--Weil--Bott theorem, a sequence $(L(p_1),\dots,L(p_n))$ is exceptional if and only if $p_j-p_i\in\Pset$ for all $i<j$, where
\[
\Pset=\{(d_1,d_2)\in\Z^2:\ 1\le d_1\le\ell-1\ \text{ or }\ 1\le d_2\le\ell-1\ \text{ or }\ d_1+d_2=\ell\}
\]
(Corollary~\ref{cor:exc-criterion}). Fullness is guaranteed, by Beilinson's theorem and Orlov's projective bundle formula, as soon as the sequence contains $\ell$ consecutive \emph{complete} columns or rows, that is, columns or rows containing $\ell+1$ consecutive lattice points (Lemma~\ref{lem:complete}).

\subsection*{Reduction to a descent.}
Call two adjacent terms of a sequence a \emph{descent} if the second coordinate drops. Rotations by the anticanonical bundle and transpositions of adjacent mutually orthogonal terms preserve both exceptionality and fullness. Using them, a non-full sequence of maximal length without non-orthogonal descents can be sorted by its second coordinate and compressed until it consists of $\ell$ complete rows, a contradiction. Thus Theorem~\ref{thm:main} reduces to the \emph{descent criterion} (Theorem~\ref{thm:descent}): a maximal-length exceptional sequence of line bundles that contains a non-orthogonal descent is full.

\subsection*{Counting.}
After rotation and translation, a non-orthogonal descent becomes the initial pair $(0,q)$, and the remaining $\ell(\ell+1)-2$ terms lie in the region $S_q=\Pset\cap(q+\Pset)$. The vector $q$ has one of two shapes (Lemma~\ref{lem:descent-shapes}), a diagonal shape and a strip shape. For all but one of them we show that $S_q$ cannot accommodate the remaining terms. The main tool is the \emph{ordered capacity} $\aord(U)$, the largest size of a subset of $U$ that can be ordered into an exceptional sequence. We bound it by partitioning the relevant part of $S_q$, which we call the \emph{wing}, into blocks of one, two, or four points. Each four-point block contains a forced three-cycle, and its fourth point is incompatible with two vertices of the cycle; hence its ordered capacity is at most two(Section~\ref{sec:nonbasic}).

\subsection*{Organization.}
Section~\ref{sec:prelim} collects the geometric input and the lattice combinatorics. Section~\ref{sec:descents} reduces Theorem~\ref{thm:main} to the descent criterion. Section~\ref{sec:nonbasic} excludes all descents except the basic ones $q=(1,-k)$. Section~\ref{sec:basic} classifies the wings of basic descents, and Section~\ref{sec:generation} proves that the resulting sequences are full. Section~\ref{sec:proof} assembles the proof.

\subsection*{Conventions.}
All varieties are over $\C$. All pullbacks, pushforwards and tensor products between derived categories are derived. For a set $S$ of objects of a triangulated category, $\gen{S}$ denotes the smallest full triangulated subcategory that contains $S$ and is closed under direct summands. All intervals $[a,b]=\{n\in\Z: a\le n\le b\}$ consist of integers; such an interval, or an index range with reversed endpoints, is empty. Throughout, $\ell\ge3$ unless stated otherwise, and
\[
h=\ell-1,\qquad N=\ell(\ell+1),\qquad I=[1,h],\qquad u=(\ell,\ell)\in\Z^2.
\]
For $d=(d_1,d_2)\in\Z^2$ we write $\sigma(d)=d_1+d_2$. We call $d_1$ the \emph{first} (horizontal) and $d_2$ the \emph{second} (vertical) coordinate. A \emph{column} is a set $\{a\}\times\Z$, a \emph{row} is a set $\Z\times\{b\}$, and an \emph{antidiagonal} is a set $\Delta_c=\{d:\sigma(d)=c\}$.

\subsection*{AI disclosure.}
I developed the mathematical ideas and proof methods in this paper, with limited AI assistance in working out specific details. In searching for four-point configurations, I observed that there was considerable freedom in choosing these configurations and that the proof did not require a sharp bound. I supplied the guiding ideas and asked GPT-6 to find a few representative configurations sufficient to obtain the required bound. After completing the research, I wrote an initial draft and used Opus 5.5 to rewrite the manuscript. I then reviewed and edited the resulting text. I take full responsibility for the mathematical content and the final manuscript.

\section{Preliminaries}\label{sec:prelim}

\subsection{The incidence variety}\label{subsec:variety}

Let $V=\C^{\ell+1}$. Let $\PP(V)$ be the projective space of lines in $V$ and $\PP(V^*)$ the projective space of hyperplanes in $V$. Set
\[
X=X_\ell=\{(\lambda,H)\in\PP(V)\times\PP(V^*):\ \lambda\subset H\},
\]
let $p_1,p_2$ be the two projections, and put
\[
L(a,b)=p_1^*\cO_{\PP(V)}(a)\otimes p_2^*\cO_{\PP(V^*)}(b),\qquad (a,b)\in\Z^2 .
\]
Then $X$ is a smooth divisor of bidegree $(1,1)$ in $\PP(V)\times\PP(V^*)$, so $\dim X=2\ell-1$ and, by adjunction,
\begin{equation}\label{eq:canonical}
\omega_X=L(-\ell,-\ell)=L(-u).
\end{equation}

Let $Y=\PP(V^*)$ and $U=H^0(Y,\cO_Y(1))\cong V$. The Euler sequence of $Y$ reads
\begin{equation}\label{eq:euler}
0\fto E\fto U\otimes\cO_Y\xrightarrow{\ \varepsilon\ }\cO_Y(1)\fto0,\qquad E=\Omega_Y(1),\quad \rk E=\ell,\quad \det E=\cO_Y(-1),
\end{equation}
and the fibre of $E$ at a point $H\in Y$ is the hyperplane $H\subset V$. With the convention that $\PP(E)$ parametrizes lines in the fibres of $E$, the second projection identifies
\begin{equation}\label{eq:projbundle}
\pi=p_2\colon X=\PP(E)\fto Y,\qquad \cO_\pi(-1)=L(-1,0),\qquad \pi_*\cO_\pi(1)=E^\vee ,
\end{equation}
where $\cO_\pi(-1)\subset\pi^*E$ is the tautological line subbundle, whose fibre at $(\lambda,H)$ is $\lambda$. Symmetrically, $p_1$ makes $X$ a $\PP^{h}$-bundle over $\PP(V)$ with relative hyperplane bundle $L(0,1)$. Consequently $\Pic(X)=\Z^2$, with $L(a,b)$ corresponding to $(a,b)$, and we freely identify line bundles with lattice points. The relative dualizing sheaf of $\pi$ is
\begin{equation}\label{eq:omega-pi}
\omega_\pi=\cO_\pi(-\ell)\otimes\pi^*\det E^\vee=L(-\ell,1).
\end{equation}

\subsection{Cohomology of line bundles}\label{subsec:cohomology}

Define
\[
\Pset=\Pset_\ell=\{d\in\Z^2:\ d_1\in I\ \text{ or }\ d_2\in I\ \text{ or }\ \sigma(d)=\ell\}.
\]

\begin{lemma}\label{lem:cohomology}
For $(a,b)\in\Z^2$, one has $H^\bullet(X,L(a,b))=0$ if and only if $(-a,-b)\in\Pset$, that is, if and only if
\[
a\in[-h,-1]\quad\text{or}\quad b\in[-h,-1]\quad\text{or}\quad a+b=-\ell .
\]
Moreover, $H^\bullet(X,\cO_X)=\C$.
\end{lemma}

\begin{proof}
This follows from \cite[Theorem~6.1 and Corollary~6.2]{AHW25}.
\end{proof}

\begin{corollary}\label{cor:exc-criterion}
Every line bundle on $X$ is exceptional. For $p_1,\dots,p_n\in\Z^2$, the sequence $(L(p_1),\dots,L(p_n))$ is exceptional if and only if
\[
p_j-p_i\in\Pset\qquad\text{for all }1\le i<j\le n .
\]
Two line bundles $L(p)$ and $L(p')$ are mutually orthogonal, that is, $\Ext^\bullet$ between them vanishes in both directions, if and only if $p'-p\in\Pset\cap(-\Pset)$.
\end{corollary}

\begin{proof}
This follows from $\Ext^\bullet(L(p_j),L(p_i))=H^\bullet(X,L(p_i-p_j))$ and Lemma~\ref{lem:cohomology}.
\end{proof}

\begin{lemma}[Duality]\label{lem:duality}
One has $\Pset=u-\Pset$; that is, $d\in\Pset$ if and only if $u-d\in\Pset$.
\end{lemma}

\begin{proof}
The three conditions defining $\Pset$ are invariant under $d\mapsto u-d$: indeed $d_i\in I$ if and only if $\ell-d_i\in I$, and $\sigma(d)=\ell$ if and only if $\sigma(u-d)=\ell$. (This is Serre duality, $H^i(L(-d))\cong H^{2\ell-1-i}(L(d-u))^\vee$, combined with~\eqref{eq:canonical}.)
\end{proof}

\subsection{Lattice combinatorics}\label{subsec:lattice}

\begin{definition}\label{def:admissible}
\begin{enumerate}[label=(\alph*)]
\item A finite sequence $\mathbf p=(p_1,\dots,p_n)$ of points of $\Z^2$ is \emph{admissible} if $p_j-p_i\in\Pset$ for all $i<j$. By Corollary~\ref{cor:exc-criterion}, this holds if and only if $(L(p_1),\dots,L(p_n))$ is exceptional. In particular, the terms of an admissible sequence are pairwise distinct, since $0\notin\Pset$. We call $\mathbf p$ \emph{full} if the corresponding exceptional sequence is full. An admissible sequence has \emph{maximal length} if $n=N$.
\item Two distinct points $p,p'$ are \emph{compatible} if $p'-p\in\Pset\cup(-\Pset)$, and \emph{orthogonal} if $p'-p\in\Pset\cap(-\Pset)$.
\item We write $p\fto p'$, and say that the arrow $p\fto p'$ is \emph{forced}, if $p'-p\in\Pset$ and $p-p'\notin\Pset$. A vector $d$ with $d\in\Pset$ and $-d\notin\Pset$ is called a \emph{forced direction}. If $p\fto p'$, then $p$ precedes $p'$ in every admissible sequence containing both points. A \emph{forced cycle} is a sequence of forced arrows $p_1\fto p_2\fto\dots\fto p_r\fto p_1$.
\item A finite set $U\subset\Z^2$ is \emph{admissibly orderable} if its elements can be arranged into an admissible sequence. We put
\begin{align*}
\aord(U)&=\max\{|U'|:\ U'\subseteq U\text{ admissibly orderable}\},\\
\acomp(U)&=\max\{|U'|:\ U'\subseteq U\text{ pairwise compatible}\}.
\end{align*}
We call $\aord(U)$ the \emph{ordered capacity} of $U$.
\end{enumerate}
\end{definition}

In the terminology of~\cite[Section~3.1]{AHW25}, admissibly orderable sets are exceptional sets of line bundles, and the forced arrows generate their exceptional poset. We shall use the following elementary facts.

\begin{lemma}\label{lem:capacity}
Let $U\subset\Z^2$ be finite.
\begin{enumerate}[label=(\alph*)]
\item $U$ is admissibly orderable if and only if its points are pairwise compatible and $U$ contains no forced cycle.
\item $\aord(U)\le\acomp(U)$, and $\aord(U)\le\sum_\nu\aord(U_\nu)$ for every partition $U=\bigsqcup_\nu U_\nu$.
\item If $p,p'$ are incompatible, then $\aord(\{p,p'\})=1$.
\item If $a\fto b\fto c\fto a$ is a forced cycle and $d$ is incompatible with two of the points $a,b,c$, then $\aord(\{a,b,c,d\})\le2$.
\end{enumerate}
\end{lemma}

\begin{proof}
(a) If $U$ is admissibly ordered, any two of its points are compatible and every forced arrow points forward, so there is no forced cycle. Conversely, suppose that the points of $U$ are pairwise compatible and that the forced arrows contain no cycle. Choose a total order on $U$ extending the forced arrows. If $p$ precedes $p'$ and $p'-p\notin\Pset$, then compatibility gives $p-p'\in\Pset$, so $p'\fto p$; this contradicts the choice of the order. Hence the order is admissible. Parts (b) and (c) are immediate. For (d), an admissibly orderable subset cannot contain all of $a,b,c$ by (a). If it omits $d$, it has at most two points. If it contains $d$, it contains neither of the two points incompatible with $d$, so again it has at most two points.
\end{proof}

\begin{lemma}\label{lem:tests}
Let $c\in\Z\sm\{0\}$ and $d\in\Z^2$.
\begin{enumerate}[label=(\alph*)]
\item The vectors $(c,0)$ and $(0,c)$ belong to $\Pset\cup(-\Pset)$ if and only if $|c|\le\ell$. For $1\le c\le\ell$, they are forced directions.
\item The vector $(c,-c)$ belongs to $\Pset\cup(-\Pset)$ if and only if $|c|\le h$, in which case it belongs to $\Pset\cap(-\Pset)$.
\item If $|d_1|>h$, $|d_2|>h$ and $\sigma(d)\notin\{\ell,-\ell\}$, then $d\notin\Pset\cup(-\Pset)$.
\item One has $-d\notin\Pset$ if and only if $d_1\notin[-h,-1]$, $d_2\notin[-h,-1]$ and $\sigma(d)\ne-\ell$. In particular, $-d\notin\Pset$ whenever $d\ne 0$ has nonnegative coordinates.
\end{enumerate}
\end{lemma}

\begin{proof}
All four statements follow directly from the definition of $\Pset$. For (a), note that $(c,0)\in\Pset$ if and only if $c\in I$ or $c=\ell$, while $(-c,0)\in\Pset$ is impossible for $c\ge1$.
\end{proof}

When checking that a displayed arrow $p\fto p'$ is forced, we verify that $d=p'-p\in\Pset$ and that the \emph{reverse vector} $-d$ has no coordinate in $I$ and coordinate sum different from $\ell$; by Lemma~\ref{lem:tests}(d) the latter is automatic when $d$ has nonnegative coordinates.

\begin{corollary}\label{cor:rows}
Let $U$ be a pairwise compatible set.
\begin{enumerate}[label=(\alph*)]
\item Every column and every row contains at most $\ell+1$ points of $U$, with equality only if these points are consecutive. The points of $U$ in a column (resp.\ row) are linearly ordered by forced arrows according to their second (resp.\ first) coordinate.
\item On an antidiagonal $\Delta_c$, the first coordinates of the points of $U\cap\Delta_c$ lie in an interval of length $h$, i.e.\ an interval containing $h+1$ integers, and any two of these points are orthogonal.
\end{enumerate}
\end{corollary}

\begin{proof}
By Lemma~\ref{lem:tests}(a), two points of a column are compatible if and only if their second coordinates differ by at most $\ell$, and then the arrow points upward. Hence the second coordinates of $U$ in a column have diameter at most $\ell$, which gives (a) for columns. Rows are analogous. Part (b) follows from Lemma~\ref{lem:tests}(b).
\end{proof}

\subsection{Generation and symmetries}\label{subsec:generation}

For $a\in\Z$ and $G\in\Db(Y)$ put
\[
\iota_a(G)=L(a,0)\otimes\pi^*G,\qquad \cB_a=\iota_a\Db(Y)\subseteq\Db(X).
\]
Each $\iota_a$ is fully faithful. By Beilinson's theorem~\cite{Bei78} and Orlov's projective bundle formula~\cite{Orl92} (see also~\cite[Chapter~8]{Huy06}), for all $r,a\in\Z$,
\begin{equation}\label{eq:SOD}
\Db(Y)=\gen{\cO_Y(r),\dots,\cO_Y(r+\ell)},\qquad \Db(X)=\gen{\cB_a,\cB_{a+1},\dots,\cB_{a+h}},
\end{equation}
where the second equality is a semiorthogonal decomposition. Applying the same results to $p_1$, with the roles of the two coordinates exchanged, gives the analogous statements for rows. In particular, $K_0(X)\cong\Z^N$ and every full exceptional sequence on $X$ has length $N$.

A column $\{a\}\times\Z$ of a set $U\subset\Z^2$ (that is, the set of points of $U$ with first coordinate $a$) is called \emph{complete} if it contains $\ell+1$ consecutive lattice points $(a,b),(a,b+1),\dots,(a,b+\ell)$. Complete rows are defined in the same way.

\begin{lemma}\label{lem:complete}
Let $U\subset\Z^2$ and let $\cC=\gen{L(p):p\in U}$.
\begin{enumerate}[label=(\alph*)]
\item If column $a$ of $U$ is complete, then $\cB_a\subseteq\cC$.
\item If $U$ has $\ell$ consecutive complete columns, or $\ell$ consecutive complete rows, then $\cC=\Db(X)$.
\end{enumerate}
\end{lemma}

\begin{proof}
If column $a$ contains $(a,b),\dots,(a,b+\ell)$, then $\cC$ contains $\iota_a(\cO_Y(b+j))$ for $0\le j\le\ell$, hence $\iota_a\Db(Y)=\cB_a$ by~\eqref{eq:SOD}. Part (b) follows from~\eqref{eq:SOD} and its analogue for $p_1$.
\end{proof}

\begin{lemma}\label{lem:operations}
Let $\mathbf p=(p_1,\dots,p_n)$ be an admissible sequence. Each of the following operations produces an admissible sequence of the same length, and the new sequence is full if and only if $\mathbf p$ is full:
\begin{enumerate}[label=(\roman*)]
\item \emph{translation} $(p_1,\dots,p_n)\mapsto(p_1+c,\dots,p_n+c)$ for $c\in\Z^2$;
\item \emph{rotation} $(p_1,\dots,p_n)\mapsto(p_2,\dots,p_n,p_1+u)$ and \emph{inverse rotation} $(p_1,\dots,p_n)\mapsto(p_n-u,p_1,\dots,p_{n-1})$;
\item \emph{dual reversal} $(p_1,\dots,p_n)\mapsto(-p_n,\dots,-p_1)$;
\item \emph{transposition} of two adjacent orthogonal terms.
\end{enumerate}
\end{lemma}

\begin{proof}
Parts (i) and (iii) follow from tensoring by a line bundle
and derived duality, respectively.
Since $\omega_X=L(-u)$, part (ii) follows from
\cite[Proposition~2.12]{AHW25}.
Part (iv) follows from
\cite[Remark~2.5 and Proposition~2.8]{AHW25}.
\end{proof}

\section{Descents and the reduction}\label{sec:descents}

\subsection{The descent criterion}

\begin{definition}\label{def:descent}
Let $\mathbf p=(p_1,\dots,p_n)$ be an admissible sequence. An adjacent pair $(p_i,p_{i+1})$ is a \emph{descent} if the second coordinate of $p_{i+1}$ is smaller than that of $p_i$. It is a \emph{non-orthogonal descent} if moreover $p_i$ and $p_{i+1}$ are not orthogonal.
\end{definition}

The heart of the paper is the following criterion, proved in Section~\ref{sec:proof}.

\begin{theorem}[Descent criterion]\label{thm:descent}
Let $\ell\ge3$. Every admissible sequence of maximal length that contains a non-orthogonal descent is full.
\end{theorem}

\begin{proof}[Proof of Theorem~\ref{thm:main}, assuming Theorem~\ref{thm:descent}]
Let $\mathbf p$ be an admissible sequence of maximal length and suppose, for a contradiction, that $\mathbf p$ is not full. By Lemma~\ref{lem:operations}, every sequence obtained from $\mathbf p$ by translations, rotations and transpositions of adjacent orthogonal terms is again admissible of maximal length and not full. By Theorem~\ref{thm:descent}, none of these sequences contains a non-orthogonal descent. For a sequence $\mathbf p'$ write $y(p'_i)$ for the second coordinate of its $i$-th term.

\emph{Sorting.} If $\mathbf p'$ has a descent, then the descent is orthogonal, and transposing it decreases the number of inversions $\#\{i<j: y(p'_i)>y(p'_j)\}$ by one. After finitely many transpositions we obtain a \emph{sorted} sequence, i.e.\ one with $y(p'_1)\le\dots\le y(p'_N)$.

\emph{Narrowing.} For a sorted sequence $\mathbf p'$ let $w(\mathbf p')=y(p'_N)-y(p'_1)$ and $\Phi(\mathbf p')=\sum_i y(p'_i)$. Let $M$ be the largest second coordinate of the first sorted sequence obtained. If $w(\mathbf p')>\ell$, rotate once and sort again. The rotated term has second coordinate $y(p'_1)+\ell<y(p'_N)\le M$, and sorting does not change the multiset of second coordinates. Hence $\Phi$ increases by $\ell$, while all second coordinates remain at most $M$. Since $\Phi\le NM$, after finitely many steps we reach a sorted sequence $\mathbf p'$ with $w(\mathbf p')\le\ell$. If $w(\mathbf p')=\ell$, let $p'_1,\dots,p'_r$ be the terms with minimal second coordinate; they form an initial segment, and $r<N$. Rotating $r$ times moves them to the end with second coordinate $y(p'_1)+\ell=y(p'_N)$. The result is sorted and has width at most $\ell-1$.

\emph{Conclusion.} We now have a sorted sequence whose terms lie in at most $\ell$ consecutive rows. By Corollary~\ref{cor:rows}, each row contains at most $\ell+1$ terms. Since $N=\ell(\ell+1)$, there are exactly $\ell$ consecutive rows, each consisting of $\ell+1$ consecutive points. By Lemma~\ref{lem:complete}, the sequence is full, a contradiction.
\end{proof}

\subsection{Shapes of non-orthogonal descents}

\begin{lemma}\label{lem:descent-shapes}
If $(p_i,p_{i+1})$ is a non-orthogonal descent of an admissible sequence, then $q=p_{i+1}-p_i$ has one of the following forms:
\begin{enumerate}[label=\textup{(\Alph*)}]
\item\label{shape:diag} $q=(\ell+k,-k)$ with $k\ge\ell$;
\item\label{shape:strip} $q=(s,-k)$ with $1\le s\le h$, $k\ge\ell$ and $k\ne\ell+s$.
\end{enumerate}
\end{lemma}

\begin{proof}
Write $q=(x,-k)$ with $k>0$. Non-orthogonality and admissibility mean $q\in\Pset$ and $-q\notin\Pset$. If $k\le h$, then $-q\in\Pset$. Hence $k\ge\ell$, and $q\in\Pset$ forces $x\in I$ or $x-k=\ell$. If $x\in I$, then $-q=(-x,k)\notin\Pset$ if and only if $k-x\ne\ell$. If $x=\ell+k$, then $-q=(-\ell-k,k)$ has no coordinate in $I$ and coordinate sum $-\ell$, so $-q\notin\Pset$.
\end{proof}

We call the descents of shape~\ref{shape:strip} with $s=1$, i.e.\ $q=(1,-k)$ with $k\ge\ell$ and $k\ne\ell+1$, \emph{basic}.

\subsection{Normalization}

Let $\mathbf p$ be an admissible sequence of maximal length with a non-orthogonal descent $(p_i,p_{i+1})$, and put $q=p_{i+1}-p_i$. Rotating $i-1$ times and translating by $-p_i$, we obtain a sequence of the form
\begin{equation}\label{eq:normalized}
\mathbf p=(0,q,p_3,\dots,p_N),\qquad p_3,\dots,p_N\in S_q:=\{p\in\Z^2:\ p\in\Pset,\ p-q\in\Pset\}.
\end{equation}
We call such a sequence \emph{normalized at $q$}. By Lemma~\ref{lem:operations}, it suffices to prove Theorem~\ref{thm:descent} for normalized sequences. The points $p_3,\dots,p_N$ are pairwise distinct and admissibly ordered.

Let $q=(s,-k)$ be of shape~\ref{shape:strip}. If $s+1\le x\le h$, then $x\in I$ and $x-s\in I$, so the whole column $\{x\}\times\Z$ lies in $S_q$. We call these columns the \emph{central columns}; there are none if $s=h$. Put
\begin{equation}\label{eq:wing}
W_q=S_q\sm\big([s+1,h]\times\Z\big),\qquad \Wing(\mathbf p)=W_q\cap\{p_3,\dots,p_N\}.
\end{equation}
We will see that $W_q$ is finite. By Corollary~\ref{cor:rows}, each central column contains at most $\ell+1$ terms of $\mathbf p$. Hence
\begin{equation}\label{eq:wing-lower}
\aord(W_q)\ \ge\ |\Wing(\mathbf p)|\ \ge\ N-2-(h-s)(\ell+1).
\end{equation}

\begin{lemma}[The duality involution]\label{lem:tau}
For $q\in\Z^2$ define $\tau_q(p)=q+u-p$.
\begin{enumerate}[label=(\alph*)]
\item $\tau_q$ is an involution of $S_q$ satisfying $\tau_q(p)-\tau_q(p')=p'-p$. It preserves compatibility and orthogonality, reverses forced arrows, maps admissibly orderable sets to admissibly orderable sets, and $\aord(\tau_q(U))=\aord(U)$ for finite $U\subseteq S_q$.
\item If $q=(s,-k)$ has shape~\ref{shape:strip}, then $\tau_q$ maps column $x$ to column $s+\ell-x$; in particular it permutes the central columns and preserves $W_q$.
\item Let $\mathbf p=(0,q,\mathbf t)$ be normalized at $q$, with tail $\mathbf t=(p_3,\dots,p_N)$. Then
\[
\mathbf p^\tau=(0,q,\tau_q(\mathbf t)^\rev),\qquad \tau_q(\mathbf t)^\rev=(\tau_q(p_N),\dots,\tau_q(p_3)),
\]
is again normalized at $q$, it is full if and only if $\mathbf p$ is full. If $q$ has shape (B), then $\Wing(\mathbf p^\tau)=\tau_q(\Wing(\mathbf p))$.
\end{enumerate}
\end{lemma}

\begin{proof}
(a) By Lemma~\ref{lem:duality}, $\tau_q(p)\in\Pset$ if and only if $p-q\in\Pset$, and $\tau_q(p)-q=u-p\in\Pset$ if and only if $p\in\Pset$. The remaining assertions follow from $\tau_q(p)-\tau_q(p')=p'-p$. (b) is immediate. (c) Put $c=q+u$. Applying dual reversal and then translating by $c$ transforms $(0,q,\mathbf t)$ into $(c-\mathbf t^\rev,\,u,\,c)$. Two inverse rotations then give
\[
(c-\mathbf t^\rev,u,c)\ \mapsto\ (q,\,c-\mathbf t^\rev,\,u)\ \mapsto\ (0,\,q,\,c-\mathbf t^\rev).
\]
Since $c-p=\tau_q(p)$, this is $\mathbf p^\tau$. Lemma~\ref{lem:operations} gives the assertions on admissibility and fullness, and (b) gives the assertion on wings.
\end{proof}

\section{Excluding non-basic descents}\label{sec:nonbasic}

\subsection{The diagonal shape}

\begin{proposition}\label{prop:diagonal}
Let $q=(\ell+k,-k)$ with $k\ge\ell$. Then $\acomp(S_q)\le\ell^2-1<N-2$. In particular, no admissible sequence of maximal length is normalized at a descent of shape~\ref{shape:diag}.
\end{proposition}

\begin{proof}
Put $a=\ell+k\ge2\ell$. A point $p=(x,y)$ lies in $\Pset$ if and only if $x\in I$, $y\in I$ or $x+y=\ell$, and $p-q\in\Pset$ if and only if $x\in a+I$, $y\in I-k$ or $x+y=2\ell$, since $\sigma(q)=\ell$. Of the nine intersections, three are empty, because $I\cap(a+I)=\emptyset$, $I\cap(I-k)=\emptyset$ and $\ell\ne2\ell$. The other six give a disjoint decomposition $S_q=R_-\sqcup R_+\sqcup\Gamma_{1,L}\sqcup\Gamma_{1,R}\sqcup\Gamma_{2,L}\sqcup\Gamma_{2,R}$, where
\begin{align*}
R_-&=I\times(I-k), & R_+&=(a+I)\times I,\\
\Gamma_{1,L}&=\{(i,2\ell-i):i\in I\}, & \Gamma_{1,R}&=\{(2\ell-j,j):j\in I\},\\
\Gamma_{2,L}&=\{(a-j,-k+j):j\in I\}, & \Gamma_{2,R}&=\{(a+i,-k-i):i\in I\}.
\end{align*}
The sets $\Gamma_{1,L}\cup\Gamma_{1,R}$ and $\Gamma_{2,L}\cup\Gamma_{2,R}$ lie on the antidiagonals $\Delta_{2\ell}$ and $\Delta_\ell$. The first coordinates of their points form the sets
\[
[1,h]\cup[\ell+1,2\ell-1]\qquad\text{and}\qquad [a-h,a-1]\cup[a+1,a+h],
\]
each consisting of two intervals of $h$ integers separated by a single missing integer. An interval containing $h+1$ integers therefore meets each of these sets in at most $h$ integers. By Corollary~\ref{cor:rows}(b), each antidiagonal contributes at most $h$ points to a compatible subset.

Let $d=(d_1,d_2)$ be the difference from a point of $R_-$ to a point of $R_+$. Then
\[
d_1\ge a+1-h=k+2>h,\qquad d_2\ge1-(h-k)=k-h+1\ge2,\qquad \sigma(d)\ge2k-\ell+4>\ell .
\]
Thus $-d\notin\Pset$, and $d\in\Pset$ if and only if $d_2\le h$. If both rectangles contain points of a compatible subset, the selected rows therefore lie in an interval of $h+1$ consecutive integers. This interval contains row $0$, which meets neither rectangle, because the rows of $R_-$ are negative and those of $R_+$ are positive. Hence at most $h$ rows are used, each containing at most $h$ points of $R_-\cup R_+$. If only one rectangle contributes, its contribution is at most $h^2$ as well. Altogether,
\[
\acomp(S_q)\le h^2+2h=\ell^2-1<\ell^2+\ell-2=N-2 .
\]
The last assertion follows because the $N-2$ terms $p_3,\dots,p_N$ of~\eqref{eq:normalized} form a pairwise compatible subset of $S_q$.
\end{proof}

\subsection{Strip shapes: the wing}\label{subsec:strip-wing}

Throughout Sections~\ref{subsec:strip-wing}--\ref{subsec:s2} we fix
\begin{equation}\label{eq:strip-q}
q=(s,-k),\qquad 2\le s\le h,\qquad k\ge\ell,
\end{equation}
and write $W=W_q$ and $\tau=\tau_q$. We do not need the condition $k\ne\ell+s$.

\begin{lemma}\label{lem:strip-points}
The set $W$ consists of the points
\begin{equation}\label{eq:ABCDEF}
\begin{aligned}
A_{ij}&=(i,\,j-k), & B_{ij}&=(h+i,\,j), & &\\
C_i&=(i,\,\ell+s-k-i), & E_i&=(h+i,\,1-i), & &\qquad 1\le i\le s,\quad 1\le j\le h,\\
D_j&=(\ell+s-k-j,\,j), & F_j&=(\ell+k-j,\,j-k), & &
\end{aligned}
\end{equation}
where a point carrying two labels is counted once.
\end{lemma}

\begin{proof}
For an interval $J$ write $V_J=J\times\Z$ and $H_J=\Z\times J$. By definition,
\[
S_q=\big(V_{I}\cup H_{I}\cup\Delta_\ell\big)\cap\big(V_{s+I}\cup H_{I-k}\cup\Delta_{\ell+s-k}\big).
\]
We have $V_I\cap V_{s+I}=V_{[s+1,h]}$, which consists of the central columns. Moreover $H_I\cap H_{I-k}=\emptyset$ since $k\ge\ell$, and $\Delta_\ell\cap\Delta_{\ell+s-k}=\emptyset$ since $k\ne s$. The other six intersections, with the central columns removed, are the six families in~\eqref{eq:ABCDEF}. The first coordinates of the $D_j$ are at most $s-1$, and those of the $F_j$ are at least $k+1>h$, so these points lie outside the central columns.
\end{proof}

In particular $W$ is finite. Since $N-2-(h-s)(\ell+1)=(h+3)h-(h-s)(h+2)$, inequality~\eqref{eq:wing-lower} becomes
\begin{equation}\label{eq:T}
\aord(W)\ \ge\ |\Wing(\mathbf p)|\ \ge\ T:=(s+1)h+2s
\end{equation}
for every sequence $\mathbf p$ of maximal length normalized at $q$. Our aim is to show $\aord(W)\le T-1$. We repeatedly use the following differences. 
\begin{equation}\label{eq:basic-diffs}
B_{i+1,j}-A_{ij}=(\ell,k),\qquad B_{1j}-D_j=(k-s+j,0),\qquad F_j-A_{sj}=(\ell+k-s-j,0).
\end{equation}
The vector $(\ell,k)$ is incompatible by Lemma~\ref{lem:tests}(c), since $k\ge\ell$.

\subsection{First estimates}

\begin{lemma}\label{lem:large-k}
If $k\ge\ell+s$, then $\acomp(W)\le T-1$.
\end{lemma}

\begin{proof}
In this range the six families in~\eqref{eq:ABCDEF} are pairwise disjoint and have no internal coincidences. Indeed, the second coordinates of the $C_i$ are at most $-1$, so $C_i\ne D_j$; moreover $C_i=A_{ij}$ would force $j=\ell+s-i\ge\ell$. The first coordinates of the $F_j$ are at least $k+1\ge h+s+2$, so $F_j\ne E_i$. All other coincidences are excluded by comparing coordinate ranges. Partition $W$ into the singletons $\{C_i\}$, $\{E_i\}$ ($1\le i\le s$) and the pairs
\begin{equation}\label{eq:pairs-large}
\{A_{ij},B_{i+1,j}\}\ (1\le i<s,\ 1\le j\le h),\qquad \{D_j,B_{1j}\},\quad\{F_j,A_{sj}\}\ (1\le j\le h).
\end{equation}
By~\eqref{eq:basic-diffs}, the pairs have differences $(\ell,k)$, $(k-s+j,0)$ and $(\ell+k-s-j,0)$. The last two have horizontal length at least $\ell+1$, so all pairs are incompatible by Lemma~\ref{lem:tests}. There are $(s-1)h+2h$ pairs and $2s$ singletons, hence $\acomp(W)\le(s+1)h+2s=T$.

Suppose a compatible subset $U$ attains this bound. Then $U$ contains every $C_i$ and every $E_i$ and exactly one point of each pair. Fix $1\le i<s$. The differences
\[
E_{i+1}-A_{i1}=(\ell,\,k-i-1),\qquad B_{i+1,h}-C_i=(\ell,\,k+i-s-1)
\]
have both coordinates at least $\ell$, so they are incompatible by Lemma~\ref{lem:tests}(c). Hence $U$ contains $B_{i+1,1}$ and $A_{ih}$. On the other hand, if the choice in the pairs $\{A_{ij},B_{i+1,j}\}$ changed between $j$ and $j+1$, then $U$ would contain one of the differences
\[
B_{i+1,j+1}-A_{ij}=(\ell,k+1),\qquad B_{i+1,j}-A_{i,j+1}=(\ell,k-1).
\]
Both are incompatible, since $k-1\ge\ell+1$. Thus the choice is constant in $j$, which contradicts the forced choices at $j=1$ and $j=h$. Hence $\acomp(W)\le T-1$.
\end{proof}

\begin{lemma}\label{lem:small-k}
Suppose $k=\ell+t$ with $0\le t\le s-1$, and put $m=s-t\in[1,s]$. Then
\[
\acomp(W)\le(s+2)h+t+2=T+(k+1-2s).
\]
In particular $\acomp(W)\le T-1$ if $k\le2s-2$.
\end{lemma}

\begin{proof}
Now $\ell+s-k=m$, so the $C_i$ and $D_j$ lie on $\Delta_m$, while the $E_i$ and $F_j$ lie on $\Delta_\ell$. Pair $A_{ij}$ with $B_{i+1,j}$ for $1\le i<s$ and $1\le j\le h$. These pairs are incompatible, and $2h$ rectangle points $A_{sj},B_{1j}$ remain. Hence the rectangles contribute at most $(s-1)h+2h=(s+1)h$ points to a compatible subset.

The first coordinates of the points on the two antidiagonals form the intervals
\[
\begin{aligned}
\Delta_m&:\ [1,s]\cup[m-h,m-1]=[m-h,s],\\
\Delta_\ell&:\ [h+1,h+s]\cup[h+t+2,2h+t+1]=[h+1,2h+t+1].
\end{aligned}
\]
These antidiagonal points are disjoint from the rectangles. On the left columns $1\le i\le s$, the second coordinate on $\Delta_m$ is $m-i\ge-t>h-k$, above the rectangle $A$. On the right columns $h+1\le x\le h+s$, the second coordinate on $\Delta_\ell$ is at most $0$, below the rectangle $B$.

Let $a\in W\cap\Delta_m$ and $b\in W\cap\Delta_\ell$ be compatible. Then $b-a=(d,\ell-m-d)$ with $d>0$ and $\ell-m\in I$. A direct check shows that $b-a\in\Pset\cup(-\Pset)$ if and only if $d\le h$ or $1\le d-\ell+m\le h$; in particular $d\le\ell-m+h$. Suppose now a compatible subset contains points of both antidiagonals. Then all its antidiagonal points have first coordinates in an interval of length at most $\ell-m+h$. That interval contains the $h-s$ integers $s+1,\dots,h$, which are not first coordinates of antidiagonal points. Since each first coordinate determines at most one antidiagonal point, the antidiagonals contribute at most
\[
(\ell-m+h+1)-(h-s)=h+t+2
\]
points. If only one antidiagonal contributes, Corollary~\ref{cor:rows}(b) gives at most $h+1\le h+t+2$ points. Adding the rectangle bound gives $(s+2)h+t+2$, and this equals $T+(k+1-2s)$.
\end{proof}

By Lemmas~\ref{lem:large-k} and~\ref{lem:small-k}, it remains to treat
\begin{equation}\label{eq:remaining}
k=\ell+t,\qquad m=s-t,\qquad 1\le m\le s\le h,\qquad h\ge s+m-2,
\end{equation}
which is equivalent to $\max(\ell,2s-1)\le k\le\ell+s-1$.

\subsection{Four-point configurations}\label{subsec:four-point}

In this subsection we assume~\eqref{eq:remaining} and $s\ge3$. The coincidences among the labels in~\eqref{eq:ABCDEF} are exactly
\begin{equation}\label{eq:coincide}
C_i=D_{m-i}\quad(1\le i<m),\qquad E_i=F_{h+2+t-i}\quad(t+2\le i\le s).
\end{equation}
Indeed, the two antidiagonals are disjoint from each other and from the rectangles, and on each antidiagonal a point is determined by its first coordinate. A list without repetitions consists of all points $A_{ij},B_{ij},D_j,F_j$ together with $C_m,\dots,C_s$ and $E_1,\dots,E_{t+1}$. Hence
\begin{equation}\label{eq:W-size}
|W|=2(s+1)h+2t+2 .
\end{equation}

\begin{lemma}\label{lem:partition-four}
The pairs
\begin{equation}\label{eq:pairs-four}
\begin{gathered}
\{A_{ij},B_{i+1,j}\}\ \ (1\le i<s,\ 1\le j\le h),\\
\{D_j,B_{1j}\}\ \ (m<j\le h),\qquad \{F_j,A_{sj}\}\ \ (1\le j\le h-m),
\end{gathered}
\end{equation}
are pairwise disjoint and incompatible. Together with the remaining points as singletons they partition $W$, and the resulting bound on $\aord(W)$ is $T+2$. Every $C_i$ and every $E_i$ is a singleton of this partition.
\end{lemma}

\begin{proof}
The rectangle points of distinct pairs are distinct, and the remaining points of the pairs lie on two different antidiagonals with distinct indices. By~\eqref{eq:basic-diffs}, the differences within the pairs are, up to sign, $(\ell,k)$, $(h+1-m+j,0)$ and $(2h+2-m-j,0)$. In the given ranges the horizontal lengths are at least $h+2=\ell+1$, so each pair is incompatible. There are $(s-1)h+2(h-m)=(s+1)h-2m$ pairs, and by~\eqref{eq:W-size} the bound equals
\[
|W|-\big((s+1)h-2m\big)=h+2t+2m+2+sh=T+2.
\]
The indices of the $D_j$ with $j\le m$ and of the $F_j$ with $j\ge h-m+1$ do not occur in~\eqref{eq:pairs-four}. By~\eqref{eq:coincide}, this shows that all $C_i$ and $E_i$ are singletons, including those carrying a second label.
\end{proof}

To gain three units, we replace three pairs and six singletons by three four-point blocks of ordered capacity two each.

\begin{lemma}\label{lem:G}
Let $1\le i<s$ and $1\le j\le m$ satisfy $m+1\le i+j\le h+1$ and $i+j\ne s$. Then the set
\[
G(i,j)=\{D_j,A_{ij},E_i,B_{i+1,j}\}
\]
has ordered capacity at most two. It consists of the pair $\{A_{ij},B_{i+1,j}\}$ of~\eqref{eq:pairs-four} and two singletons.
\end{lemma}

\begin{proof}
Put $\delta=i+j-m\in I$. The successive differences along $D_j\fto A_{ij}\fto E_i\fto D_j$ are
\begin{equation}\label{eq:G-diffs}
(\delta,-k),\qquad (h,\,h+2+t-i-j),\qquad (-h-\delta,\,i+j-1).
\end{equation}
Each lies in $\Pset$: the first coordinates of the first two vectors lie in $I$, and the second coordinate of the third does. The first reverse vector $(-\delta,k)$ has no coordinate in $I$ and coordinate sum $\ell+s-i-j\ne\ell$. The second reverse vector has both coordinates negative, since $i+j\le h+1$. The third reverse vector $(h+\delta,1-i-j)$ has no coordinate in $I$ and coordinate sum $\ell-m\ne\ell$. Hence all three arrows are forced. Moreover, $B_{i+1,j}-A_{ij}=(\ell,k)$ and $B_{i+1,j}-D_j=(\ell+\delta,0)$ are incompatible, so Lemma~\ref{lem:capacity}(d) applies. The last statement follows from Lemma~\ref{lem:partition-four}, since $D_j$ with $j\le m$ is either $D_m$ or a point $C_{m-j}$.
\end{proof}

The involution $\tau$ of Lemma~\ref{lem:tau} is $\tau(p)=(s+\ell,-t)-p$. On the labels it acts by
\begin{equation}\label{eq:tau-labels}
\tau(A_{ij})=B_{s+1-i,\,h+1-j},\quad \tau(D_j)=F_{h+1-j},\quad \tau(E_i)=C_{s+1-i},
\end{equation}
and, $\tau$ being an involution, $\tau(B_{ij})=A_{s+1-i,h+1-j}$, $\tau(F_j)=D_{h+1-j}$ and $\tau(C_i)=E_{s+1-i}$. In particular, $\tau$ preserves the partition of Lemma~\ref{lem:partition-four}, and $\aord(\tau(G))=\aord(G)$ for every $G\subseteq W$.

\begin{proposition}\label{prop:four-point}
Assume~\eqref{eq:remaining} and $s\ge3$. Then $\aord(W)\le T-1$.
\end{proposition}

\begin{proof}
\emph{Case $m\ge2$.} Consider the three groups
\begin{equation}\label{eq:groups-m2}
G(s-1,m),\qquad \tau\big(G(s-1,m)\big),\qquad H=\{D_{m-1},A_{s1},E_s,F_1\}.
\end{equation}
The indices $(i,j)=(s-1,m)$ satisfy the hypotheses of Lemma~\ref{lem:G}: $i+j=s+m-1$ lies in $[m+1,h+1]$ by~\eqref{eq:remaining} and $s\ge2$, and differs from $s$ because $m\ne1$. Hence the first two groups have ordered capacity at most two. In $H$, consider $D_{m-1}\fto A_{s1}\fto E_s\fto D_{m-1}$, with successive differences
\[
(s-1,\,1-h-s),\qquad (h,\,h+1-m),\qquad (1-h-s,\,m+s-2).
\]
They lie in $\Pset$ because $s-1\in I$, $h\in I$ and $m+s-2\in I$, respectively. The first and third reverse vectors have no coordinate in $I$ and coordinate sums $h$ and $\ell-m$. The second reverse vector has both coordinates negative. So the cycle is forced. The fourth point gives the differences
\[
F_1-A_{s1}=(2h+1-m,\,0),\qquad F_1-D_{m-1}=(2h+t,\,1-h-s).
\]
The first is incompatible because $2h+1-m\ge\ell+1$, which follows from $h\ge s+m-2\ge m+1$. The second is incompatible by Lemma~\ref{lem:tests}(c), since its coordinate sum is $\ell-m$. Thus $\aord(H)\le2$ by Lemma~\ref{lem:capacity}(d).

The three groups contain the pairs
\[
\{A_{s-1,m},B_{sm}\},\qquad \{A_{1,h+1-m},B_{2,h+1-m}\},\qquad \{A_{s1},F_1\}
\]
of~\eqref{eq:pairs-four}, where $1\le h-m$ for the last one, and the following six singletons:
\[
\begin{array}{c|c|c}
 & \text{points} & \text{first coordinates}\\\hline
\Delta_m & D_m,\ D_{m-1}=C_1,\ C_2 & 0,\ 1,\ 2\\
\Delta_\ell & E_{s-1},\ E_s,\ F_{h+1-m} & h+s-1,\ h+s,\ h+s+1
\end{array}
\]
The pairs are distinct because their $A$-columns $s-1$, $1$, $s$ are distinct, and the table shows that the six singletons are distinct. Hence the three groups are pairwise disjoint.

\emph{Case $m=1$.} Then $t=s-1$, and there are no coincidences among the labels. Consider
\begin{equation}\label{eq:groups-m1}
G(s-2,1),\qquad \tau\big(G(s-2,1)\big),\qquad H'=\{C_2,B_{1h},E_s,D_h\}.
\end{equation}
The indices $(s-2,1)$ satisfy the hypotheses of Lemma~\ref{lem:G}, since $s\ge3$ and $s-1\le h+1$. In $H'$, the cycle $C_2\fto B_{1h}\fto E_s\fto C_2$ has successive differences
\[
(h-1,\,h+1),\qquad (s-1,\,1-s-h),\qquad (2-h-s,\,s-2).
\]
They lie in $\Pset$ because $h-1$, $s-1$ and $s-2$ lie in $I$; note that $h\ge s-1\ge2$. The first reverse vector is negative in both coordinates. The other two reverse vectors have no coordinate in $I$ and coordinate sum $h$. So the cycle is forced. In addition,
\[
C_2-D_h=(h+1,\,-h-1),\qquad E_s-D_h=(2h+s-1,\,1-s-h)
\]
are incompatible by Lemma~\ref{lem:tests}(c), with coordinate sums $0$ and $h$. Hence $\aord(H')\le2$. The pairs used in~\eqref{eq:groups-m1} are
\[
\{A_{s-2,1},B_{s-1,1}\},\qquad \{A_{2h},B_{3h}\},\qquad \{D_h,B_{1h}\},
\]
and the six singletons are
\[
\begin{array}{c|c|c}
 & \text{points} & \text{first coordinates}\\\hline
\Delta_1 & D_1,\ C_3,\ C_2 & 0,\ 3,\ 2\\
\Delta_\ell & E_{s-2},\ F_h,\ E_s & h+s-2,\ h+s+1,\ h+s
\end{array}
\]
The first two pairs share columns only when $s=4$, and then they lie in the distinct rows $1$ and $h$. The third pair is the only one containing a $D$-point, and its $B$-point lies in column $h+1$, whereas the $B$-points of the first two pairs lie in columns $h+s-1$ and $h+3$. The table separates the six singletons, so the three groups are again pairwise disjoint.

In both cases, three pairs and six singletons of the partition of Lemma~\ref{lem:partition-four}, with total bound $9$, have been replaced by three groups of total ordered capacity at most $6$. By Lemma~\ref{lem:capacity}(b), $\aord(W)\le T+2-3=T-1$.
\end{proof}

\subsection{The case \texorpdfstring{$s=2$}{s=2}}\label{subsec:s2}

\begin{proposition}\label{prop:s2}
If $s=2$ and $k\ge\ell$, then $\aord(W)\le3h+3=T-1$.
\end{proposition}

\begin{proof}
Lemma~\ref{lem:large-k} covers $k\ge\ell+2$, so it remains to treat $k=\ell$ and $k=\ell+1$. Recall from Lemma~\ref{lem:tests} the two incompatibility tests
\begin{equation}\label{eq:tests}
\begin{aligned}
&|d|>\ell &&\Longrightarrow\ \ (d,0)\ \text{is incompatible},\\
&|d_1|,|d_2|>h\ \text{ and }\ d_1+d_2\notin\{\ell,-\ell\} &&\Longrightarrow\ \ (d_1,d_2)\ \text{is incompatible}.
\end{aligned}
\end{equation}
In all cyclic configurations below, every displayed successive difference lies in $\Pset$. Each reverse vector either has both coordinates negative, or has no coordinate in $I$ and coordinate sum in $\{h-1,h\}$; hence every displayed cycle is forced. The fourth point of each four-point group is incompatible with two cycle vertices, so each group has ordered capacity at most two by Lemma~\ref{lem:capacity}(d).

\emph{Case $k=\ell$, $h\ge3$.} Here $t=0$ and $m=2$, the only coincidences are $C_1=D_1$ and $E_2=F_h$, and $|W|=6h+2$. Consider
\begin{align*}
G_1&=\{A_{12},E_1,D_2,B_{22}\},\\
G_2&=\{A_{21},E_2,D_1,F_1\},\\
G_3&=\{C_2,B_{1h},F_{h-1},D_h\},
\end{align*}
where in each group the first three points form a cycle in the displayed order, with successive differences
\[
\begin{array}{c|ccc}
G_1,G_2 & (h,h-1) & (-h-1,2) & (1,-h-1)\\
G_3 & (h-1,h) & (2,-h-2) & (-h-1,2).
\end{array}
\]
The fourth points give the incompatible differences
\begin{align*}
B_{22}-A_{12}&=(h+1,h+1), & B_{22}-D_2&=(h+2,0),\\
F_1-A_{21}&=(2h-1,0), & F_1-D_1&=(2h,-h-1),\\
B_{1h}-D_h&=(2h-1,0), & F_{h-1}-D_h&=(2h+1,-h-2),
\end{align*}
where $2h-1>\ell$ uses $h\ge3$. The remaining points are partitioned into the pairs below, each incompatible by~\eqref{eq:tests}, and the four singletons $B_{11}$, $B_{12}$, $A_{2,h-1}$, $A_{2h}$:
\[
\begin{array}{c|c|c}
\text{pair} & \text{range} & \text{difference (up to sign)}\\\hline
\{A_{1j},B_{2j}\} & 1\le j\le h,\ j\ne2 & (h+1,h+1)\\
\{D_j,B_{1j}\} & 3\le j\le h-1 & (h-1+j,0)\\
\{F_j,A_{2j}\} & 2\le j\le h-2 & (2h-j,0)
\end{array}
\]
The three groups use the $D$-indices $1,2,h$ and the $F$-indices $1,h-1,h$ (recall $E_2=F_h$); all other $D$- and $F$-indices occur in the pairs. These indices are distinct also when $h=3$, so groups, pairs and singletons form a partition of $W$. There are $(h-1)+(h-3)+(h-3)=3h-7$ pairs, and therefore
\[
\aord(W)\le 3\cdot2+(3h-7)+4=3h+3 .
\]

\emph{Case $k=\ell+1$, $h\ge3$.} Here $t=1$ and $m=1$, all labels are distinct, and $|W|=6h+4$. Consider
\begin{align*}
H_1&=\{D_1,A_{22},E_1,F_2\},\\
H_2&=\{C_1,B_{1h},E_2,D_h\},\\
H_3&=\{C_2,B_{1,h-1},F_h,D_{h-1}\},
\end{align*}
with cycle differences
\[
\begin{array}{c|ccc}
H_1 & (2,-h-1) & (h-1,h) & (-h-1,1)\\
H_2 & (h,h) & (1,-h-1) & (-h-1,1)\\
H_3 & (h-1,h) & (2,-h-1) & (-h-1,1)
\end{array}
\]
and fourth-point differences
\begin{align*}
F_2-D_1&=(2h+1,-h-1), & F_2-A_{22}&=(2h-1,0),\\
B_{1h}-D_h&=(2h,0), & E_2-D_h&=(2h+1,-h-1),\\
B_{1,h-1}-D_{h-1}&=(2h-1,0), & F_h-D_{h-1}&=(2h+1,-h-1).
\end{align*}
The remaining points form the incompatible pairs
\[
\begin{array}{c|c|c}
\text{pair} & \text{range} & \text{difference (up to sign)}\\\hline
\{A_{1j},B_{2j}\} & 1\le j\le h & (h+1,h+2)\\
\{D_j,B_{1j}\} & 2\le j\le h-2 & (h+j,0)\\
\{F_j,A_{2j}\} & 1\le j\le h-1,\ j\ne2 & (2h+1-j,0)
\end{array}
\]
and the singletons $B_{11}$, $A_{2h}$. The groups use the $D$-indices $1,h-1,h$ and the $F$-indices $2,h$, which are distinct also when $h=3$. There are $h+(h-3)+(h-2)=3h-5$ pairs, so
\[
\aord(W)\le 3\cdot2+(3h-5)+2=3h+3 .
\]

\emph{Case $h=2$.} Then $\ell=3$, $T=10$, and $k\in\{3,4\}$. For $k=3$ the set $W$ has $14$ points, which we partition into the three forced cycles
\[
\begin{gathered}
(0,2)\fto(2,-2)\fto(3,0)\fto(0,2),\\
(1,1)\fto(3,2)\fto(4,-1)\fto(1,1),\\
(2,0)\fto(3,1)\fto(5,-2)\fto(2,0),
\end{gathered}
\]
the two pairs $\{(1,-2),(4,1)\}$ and $\{(1,-1),(4,2)\}$, and the singleton $\{(2,-1)\}$. The reverse cycle vectors having a positive coordinate are $(-2,4)$, $(3,-2)$, $(-1,3)$ and $(-2,3)$, none of which lies in $\Pset_3$; the other reverse vectors have both coordinates negative. Both pairs have difference $(3,3)$ and are incompatible. Thus $\aord(W)\le6+2+1=9$.

For $k=4$ the set $W$ has $16$ points, which we partition into the three forced cycles
\[
\begin{gathered}
(0,1)\fto(2,-2)\fto(3,0)\fto(0,1),\\
(1,0)\fto(3,2)\fto(4,-1)\fto(1,0),\\
(2,-1)\fto(3,1)\fto(5,-2)\fto(2,-1),
\end{gathered}
\]
the triple $\{(-1,2),(2,-3),(6,-3)\}$ and the two pairs $\{(1,-3),(4,1)\}$ and $\{(1,-2),(4,2)\}$. The reverse cycle vectors having a positive coordinate are $(-2,3)$, $(-1,3)$ and $(3,-1)$, none of which lies in $\Pset_3$. The triple has pairwise differences $(3,-5)$, $(7,-5)$ and $(4,0)$, all incompatible, and both pairs have difference $(3,4)$. Thus $\aord(W)\le6+1+2=9$.
\end{proof}

\begin{corollary}\label{cor:no-strip}
No admissible sequence of maximal length is normalized at $q=(s,-k)$ with $2\le s\le h$ and $k\ge\ell$.
\end{corollary}

\begin{proof}
Lemmas~\ref{lem:large-k} and~\ref{lem:small-k} and Propositions~\ref{prop:four-point} and~\ref{prop:s2} cover every parameter in~\eqref{eq:strip-q}, and give $\aord(W)\le T-1$. This contradicts~\eqref{eq:T}.
\end{proof}

\section{Basic descents: the wings}\label{sec:basic}

Throughout this section $q=(1,-k)$ with $k\ge\ell$ and $k\ne\ell+1$, and $W=W_q$. The central columns are $2\le x\le h$, and $\tau=\tau_q$ is given by $\tau(p)=(\ell+1,\ell-k)-p$.

\begin{lemma}\label{lem:basic-points}
The set $W$ is the disjoint union of the four families
\begin{equation}\label{eq:LRDF}
\begin{aligned}
L_i&=(1,\,-k+i), &&1\le i\le\ell, & R_t&=(\ell,\,t), &&0\le t\le h,\\
D_z&=(\ell+1-k-z,\,z), &&1\le z\le h, & F_j&=(\ell+k-j,\,-k+j), &&1\le j\le h.
\end{aligned}
\end{equation}
In particular $|W|=4\ell-2$. The involution $\tau$ satisfies $\tau(L_i)=R_{\ell-i}$ and $\tau(D_z)=F_{\ell-z}$.
\end{lemma}

\begin{proof}
This is the computation of Lemma~\ref{lem:strip-points} with $s=1$. Here $p-q\in\Pset$ if and only if $x\in[2,\ell]$, $y\in I-k$ or $x+y=\ell+1-k$. Outside the central columns, the condition $x\in I$ leaves column $1$, which gives the $L_i$; the condition $x\in[2,\ell]$ leaves column $\ell$, which gives the $R_t$. The intersections $H_I\cap\Delta_{\ell+1-k}$ and $\Delta_\ell\cap H_{I-k}$ give the $D_z$ and $F_j$, whose first coordinates are at most $0$ and at least $\ell+1$, respectively. The formulas for $\tau$ are immediate.
\end{proof}

\begin{lemma}\label{lem:basic-count}
One has $\aord(W)\le2\ell$. If $\mathbf p=(0,q,p_3,\dots,p_N)$ is admissible and normalized at $q$, then $|\Wing(\mathbf p)|=2\ell$, and every central column of $\mathbf p$ is complete.
\end{lemma}

\begin{proof}
First let $k>\ell$. Partition $W$ into
\begin{equation}\label{eq:basic-partition-large}
\{D_z,R_z\}\ (1\le z\le h),\qquad \{F_j,L_j\}\ (1\le j\le h),\qquad \{L_\ell\},\qquad \{R_0\}.
\end{equation}
The pairs have differences $R_z-D_z=(k+z-1,0)$ and $F_j-L_j=(\ell+k-j-1,0)$, with horizontal lengths larger than $\ell$; hence they are incompatible. This gives $\aord(W)\le\acomp(W)\le2h+2=2\ell$.

Now let $k=\ell$ and consider
\begin{equation}\label{eq:G1G2}
G_1=\{D_1,F_1,L_1,R_0\},\qquad G_2=\tau(G_1)=\{F_h,D_h,R_h,L_\ell\}.
\end{equation}
The cycle $D_1\fto L_1\fto R_0\fto D_1$ has successive differences $(1,-\ell)$, $(h,h)$ and $(-\ell,1)$. Each lies in $\Pset$. The reverse vectors $(-1,\ell)$ and $(\ell,-1)$ have no coordinate in $I$ and coordinate sum $h\ne\ell$, while $(-h,-h)$ is negative; so the cycle is forced. Moreover, $F_1$ is incompatible with $D_1$ and $L_1$, because
\[
F_1-D_1=(2\ell-1,\,-\ell),\qquad F_1-L_1=(2\ell-2,\,0),
\]
where the second vector is incompatible because $\ell\ge3$. By Lemma~\ref{lem:capacity}(d), $\aord(G_1)\le2$, and by Lemma~\ref{lem:tau}, $\aord(G_2)\le2$. Complete the partition of $W$ by
\begin{equation}\label{eq:basic-partition-l}
\{D_r,R_r\},\ \{F_r,L_r\}\ (2\le r\le h-1),\qquad \{R_1\},\qquad \{L_h\}.
\end{equation}
These pairs have horizontal distances $h+r$ and $2\ell-r-1$, both larger than $\ell$. Hence
\[
\aord(W)\le4+2(h-2)+2=2\ell .
\]
(When $h=2$ the ranges in~\eqref{eq:basic-partition-l} are empty and the same count applies.)

For the second assertion, let $n_x$ be the number of terms of $\mathbf p$ in the central column $x$. By Corollary~\ref{cor:rows} and the first part,
\[
N-2=|\Wing(\mathbf p)|+\sum_{x=2}^{h}n_x\le2\ell+(\ell-2)(\ell+1)=N-2 .
\]
Hence all inequalities are equalities: the wing has $2\ell$ points, and each central column contains $\ell+1$ terms, which are consecutive by Corollary~\ref{cor:rows}.
\end{proof}

Put
\[
T_0=\{L_1,\dots,L_\ell,R_0,\dots,R_h\}.
\]

\begin{proposition}\label{prop:basic-large}
Let $k\ge\ell+2$. Every admissible sequence $\mathbf p$ of maximal length normalized at $q=(1,-k)$ satisfies $\Wing(\mathbf p)=T_0$ and is full.
\end{proposition}

\begin{proof}
By Lemma~\ref{lem:basic-count}, $\Wing(\mathbf p)$ is an admissibly orderable subset of $W$ of size $2\ell$. Hence equality holds in the partition~\eqref{eq:basic-partition-large}: $L_\ell,R_0\in\Wing(\mathbf p)$, and each pair contributes exactly one point. Since
\[
F_1-R_0=(k-1,\,1-k)
\]
is incompatible by Lemma~\ref{lem:tests}(b), $F_1$ is absent, and so $L_1\in\Wing(\mathbf p)$. Now let $1\le z\le h$. The differences
\[
R_0-L_1=(h,\,k-1),\qquad D_z-R_0=(1-k-z,\,z)
\]
give forced arrows $L_1\fto R_0\fto D_z$. Indeed, both vectors lie in $\Pset$; the first reverse vector is negative; and the second reverse vector $(k+z-1,-z)$ has no coordinate in $I$ and coordinate sum $k-1\ne\ell$, because $k\ne\ell+1$. On the other hand,
\[
D_z-L_1=(\ell-k-z,\,k+z-1)\notin\Pset,
\]
since its first coordinate is negative, its second exceeds $h$, and its coordinate sum is $\ell-1$. So if $D_z$ and $L_1$ were compatible, we would have $D_z\fto L_1$ and hence the forced cycle $D_z\fto L_1\fto R_0\fto D_z$. Therefore no $D_z$ belongs to $\Wing(\mathbf p)$.

Applying this to the sequence $\mathbf p^\tau$ of Lemma~\ref{lem:tau}(c), whose wing is $\tau(\Wing(\mathbf p))$, we find that $\tau(\Wing(\mathbf p))$ contains no $D_z$. By Lemma~\ref{lem:basic-points}, $\Wing(\mathbf p)$ contains no $F_j$. The pairs in~\eqref{eq:basic-partition-large} then force $\Wing(\mathbf p)=T_0$.

The point $q=(1,-k)$ and $L_1,\dots,L_\ell$ form the complete column $1$, namely $(1,-k),\dots,(1,\ell-k)$. Rotate $\mathbf p$ once; this replaces the initial term $0$ by $u=(\ell,\ell)$, which together with $R_0,\dots,R_h$ completes column $\ell$. The central columns $2,\dots,h$ are complete by Lemma~\ref{lem:basic-count}. Thus the rotated sequence has the $\ell$ consecutive complete columns $1,\dots,\ell$ and is full by Lemma~\ref{lem:complete}; hence $\mathbf p$ is full by Lemma~\ref{lem:operations}.
\end{proof}

\subsection{The boundary value \texorpdfstring{$k=\ell$}{k=l}}

For the rest of this section we take $k=\ell$, so that $q=(1,-\ell)$ and
\begin{equation}\label{eq:LRDF-l}
L_i=(1,\,i-\ell),\qquad R_t=(\ell,\,t),\qquad D_z=(1-z,\,z),\qquad F_j=(2\ell-j,\,j-\ell),
\end{equation}
and $\rho:=\tau$ is the reflection $\rho(p)=(\ell+1,0)-p$. Thus
\begin{equation}\label{eq:rho-labels}
\rho(L_i)=R_{\ell-i},\qquad \rho(R_t)=L_{\ell-t},\qquad \rho(D_z)=F_{\ell-z},\qquad \rho(F_j)=D_{\ell-j}.
\end{equation}
A direct substitution into the definition of $\Pset$ gives Table~\ref{tab:directions}. For a pair $(A,B)$ of points from two of the four families, it records when $B-A\in\Pset$ (``$A$ before $B$ is allowed'') and when $A-B\in\Pset$ (``$B$ before $A$ is allowed''). Within the $L$-family and within the $R$-family the indices are forced to increase, by Corollary~\ref{cor:rows}(a). The $D$-points are pairwise orthogonal, and so are the $F$-points, by Lemma~\ref{lem:tests}(b).

\begin{table}[ht]
\centering
\renewcommand{\arraystretch}{1.15}
\begin{tabular}{llll}
\toprule
$(A,B)$ & $B-A$ & $B-A\in\Pset$ & $A-B\in\Pset$\\
\midrule
$(L_i,R_t)$ & $(h,\,t+\ell-i)$ & always & never\\
$(D_z,L_i)$ & $(z,\,i-\ell-z)$ & always & $i>z$\\
$(R_t,F_j)$ & $(\ell-j,\,j-\ell-t)$ & always & $t<j$\\
$(D_z,R_t)$ & $(h+z,\,t-z)$ & $t>z$ or $t=1$ & $z>t$\\
$(L_i,F_j)$ & $(2\ell-j-1,\,j-i)$ & $j>i$ or $i=h$ & $i>j$\\
$(D_z,F_j)$ & $(2\ell-1+z-j,\,j-\ell-z)$ & never & $j>z$\\
\bottomrule
\end{tabular}
\vspace{6pt}
\caption{Allowed relative orders between the four families for $k=\ell$.}
\label{tab:directions}
\end{table}

For $1\le m\le h$ define
\begin{equation}\label{eq:Gamma}
\Gamma_m=\{L_1,\dots,L_\ell\}\cup\{D_1,\dots,D_m\}\cup\{R_1,R_{m+1},\dots,R_h\},
\end{equation}
and define
\begin{equation}\label{eq:T-sets}
\begin{aligned}
T_{00}&=T_0, & T_{10}&=(T_0\sm\{L_1\})\cup\{D_1\},\\
T_{01}&=(T_0\sm\{R_h\})\cup\{F_h\}, & T_{11}&=(T_0\sm\{L_1,R_h\})\cup\{D_1,F_h\}.
\end{aligned}
\end{equation}
Note that $\rho(T_{00})=T_{00}$, $\rho(T_{11})=T_{11}$ and $\rho(T_{10})=T_{01}$.

\begin{theorem}\label{thm:classification}
Let $k=\ell$. The admissibly orderable subsets of $W$ of cardinality $2\ell$ are exactly
\[
\Gamma_m,\quad \rho(\Gamma_m)\quad(1\le m\le h),\qquad T_{00},\quad T_{10},\quad T_{01},\quad T_{11}.
\]
These $2\ell+2$ sets are pairwise distinct.
\end{theorem}

\begin{proof}
Let $T\subseteq W$ be admissibly orderable with $|T|=2\ell$. Then the partition of $W$ into $G_1$, $G_2$ and the blocks~\eqref{eq:basic-partition-l} used in the proof of Lemma~\ref{lem:basic-count} is saturated by $T$. That is, $R_1,L_h\in T$, each pair in~\eqref{eq:basic-partition-l} contributes exactly one point, and each of $G_1,G_2$ contributes exactly two points. The admissibly orderable two-element subsets of $G_1$ are
\begin{equation}\label{eq:G1-choices}
\{D_1,L_1\},\qquad \{D_1,R_0\},\qquad \{L_1,R_0\},\qquad \{F_1,R_0\},
\end{equation}
because $F_1$ is incompatible with $D_1$ and $L_1$. By~\eqref{eq:rho-labels}, the admissibly orderable two-element subsets of $G_2$ are their images $\{F_h,R_h\}$, $\{F_h,L_\ell\}$, $\{R_h,L_\ell\}$, $\{D_h,L_\ell\}$.

\emph{Step 1: $T$ contains both a $D$-point and an $F$-point.} Suppose $D_z,F_j\in T$. Table~\ref{tab:directions} requires $j>z$. If $z\ge2$, then $j\ge3$, so $h\ge3$. By the table, $F_2$ is incompatible with $D_z$, so the saturated pair $\{F_2,L_2\}$ gives $L_2\in T$. The table then yields the forced cycle
\[
D_z\fto L_2\fto F_j\fto D_z,
\]
which is impossible. Thus $D_1$ is the only $D$-point that can occur together with an $F$-point. Applying this to $\rho(T)$, we see that $F_h=\rho(D_1)$ is the only $F$-point that can occur together with a $D$-point. Hence $D_1,F_h\in T$, and $T$ contains no other $D$- or $F$-points. The forced cycles
\[
D_1\fto L_1\fto F_h\fto D_1,\qquad F_h\fto D_1\fto R_h\fto F_h
\]
exclude $L_1$ and $R_h$; both cycles use $h\ge2$. Saturation now gives $G_1\cap T=\{D_1,R_0\}$, $G_2\cap T=\{F_h,L_\ell\}$, $L_r,R_r\in T$ for $2\le r\le h-1$, and $R_1,L_h\in T$. Hence $T=T_{11}$.

\emph{Step 2: $T$ contains a $D$-point but no $F$-point.} The saturated pairs $\{F_r,L_r\}$ and the singleton $\{L_h\}$ give $L_2,\dots,L_h\in T$. Since $F_h\notin T$, the choice in $G_2$ is $\{R_h,L_\ell\}$ or $\{D_h,L_\ell\}$, so $L_\ell\in T$. If $D_1\notin T$, then~\eqref{eq:G1-choices} forces $\{L_1,R_0\}\subseteq T$, and any $D_z\in T$ would give the forced cycle
\[
D_z\fto L_1\fto R_0\fto D_z .
\]
Hence $D_1\in T$. Suppose first $G_1\cap T=\{D_1,R_0\}$. For $z\ge2$, the forced cycle
\[
D_z\fto L_2\fto R_0\fto D_z
\]
excludes $D_z$. Saturation then gives $R_2,\dots,R_{h-1}\in T$ and $G_2\cap T=\{R_h,L_\ell\}$, so $T=T_{10}$. Otherwise $G_1\cap T=\{D_1,L_1\}$, and $T$ contains every $L_i$ but not $R_0$. Let $m=\max\{z:D_z\in T\}$. For $2\le t<m$, the forced cycle
\[
D_m\fto L_1\fto R_t\fto D_m
\]
excludes $R_t$, and the saturated pair $\{D_t,R_t\}$ gives $D_t\in T$. For $m<t\le h-1$, the maximality of $m$ and saturation give $R_t\in T$. Finally, $G_2\cap T=\{D_h,L_\ell\}$ if $m=h$, and $G_2\cap T=\{R_h,L_\ell\}$ if $m<h$. Together with $R_1\in T$, this gives $T=\Gamma_m$.

\emph{Step 3: the remaining cases.} If $T$ contains an $F$-point but no $D$-point, then $\rho(T)$ contains a $D$-point but no $F$-point. By Step~2, $T\in\{\rho(\Gamma_m):1\le m\le h\}\cup\{T_{01}\}$. If $T$ contains neither, saturation gives $T=T_{00}$.

\emph{Converse.} The following orders are admissible:
\begin{equation}\label{eq:orders}
\begin{array}{ll}
\Gamma_m: & D_1,\dots,D_m,\ L_1,\dots,L_\ell,\ R_1,R_{m+1},\dots,R_h,\\
\rho(\Gamma_m): & L_1,\dots,L_{h-m},\ L_h,\ R_0,\dots,R_h,\ F_{\ell-m},\dots,F_h,\\
T_{00}: & L_1,\dots,L_\ell,\ R_0,\dots,R_h,\\
T_{10}: & L_2,\dots,L_\ell,\ R_0,\ D_1,\ R_1,\dots,R_h,\\
T_{01}: & L_1,\dots,L_h,\ F_h,\ L_\ell,\ R_0,\dots,R_{h-1},\\
T_{11}: & L_2,\dots,L_h,\ F_h,\ L_\ell,\ R_0,\ D_1,\ R_1,\dots,R_{h-1}.
\end{array}
\end{equation}
For $\Gamma_m$, $T_{00}$, $T_{10}$ and $T_{11}$, this follows directly from Table~\ref{tab:directions} and the forced increase of indices within the $L$- and $R$-families. The orders for $\rho(\Gamma_m)$ and $T_{01}$ are obtained by applying $\rho$ to the orders for $\Gamma_m$ and $T_{10}$ and reversing them; this preserves admissibility by Lemma~\ref{lem:tau}(a).

Finally, each $\Gamma_m$ contains every $L_i$, no $F_j$ and not $R_0$, and the largest $D$-index distinguishes the different values of $m$. Dually, each $\rho(\Gamma_m)$ contains every $R_t$, no $D_z$ and not $L_\ell$. Hence these $2h$ sets are pairwise distinct. The four sets in~\eqref{eq:T-sets} all contain $R_0$ and $L_\ell$, so none of them belongs to either family, and they are distinguished by the presence of $D_1$ and $F_h$. The total number is $2h+4=2\ell+2$.
\end{proof}

\section{Generation from the saturated configurations}\label{sec:generation}

Throughout this section $q=(1,-\ell)$, and $L_i,R_t,D_z,F_j$ are as in~\eqref{eq:LRDF-l}. Put $\cD=\Db(Y)$. Recall that $\iota_a(G)=L(a,0)\otimes\pi^*G$ and $\cB_a=\iota_a\cD$, and define
\[
Q\colon\Db(X)\fto\cD,\qquad Q(F)=R\pi_*\big(F\otimes L(-\ell,0)\big).
\]

\begin{lemma}\label{lem:Q}
Let $\cC\subseteq\Db(X)$ be a full triangulated subcategory closed under direct summands and containing $\cB_1,\dots,\cB_h$. If $\gen{Q(F_\nu)}_\nu=\cD$ for some objects $F_\nu\in\cC$, then $\cC=\Db(X)$.
\end{lemma}

\begin{proof}
By adjunction, $\iota_\ell$ is left adjoint to $Q$, and $Q\iota_\ell\simeq\id_\cD$ because $R\pi_*\cO_X=\cO_Y$. For $1\le a\le h$, one has $Q(\cB_a)=0$ since $R\pi_*\cO_\pi(a-\ell)=0$ for $a-\ell\in[-h,-1]$. Since $\iota_\ell$ is left adjoint to $Q$, the kernel of $Q$ is the right orthogonal $\cB_\ell^{\perp}=\{F:\Hom(\cB_\ell,F[i])=0\text{ for all }i\}$. The semiorthogonal decomposition $\Db(X)=\gen{\cB_1,\dots,\cB_h,\cB_\ell}$ from~\eqref{eq:SOD} therefore gives
\[
\ker Q=\cB_\ell^{\perp}=\gen{\cB_1,\dots,\cB_h}\subseteq\cC .
\]
For every $F\in\Db(X)$, the counit gives a triangle
\begin{equation}\label{eq:counit}
\iota_\ell Q(F)\fto F\fto A_F\fto\iota_\ell Q(F)[1],\qquad A_F\in\ker Q,
\end{equation}
since $Q$ applied to the counit is an isomorphism. Hence $\iota_\ell Q(F_\nu)\in\cC$ for all $\nu$, and therefore $\cB_\ell=\iota_\ell\cD=\gen{\iota_\ell Q(F_\nu)}_\nu\subseteq\cC$. By~\eqref{eq:SOD}, $\cC=\Db(X)$.
\end{proof}

We shall compute $Q$ on the relevant line bundles using relative Serre duality for $\pi$ and~\eqref{eq:omega-pi}:
\begin{equation}\label{eq:rel-duality}
\begin{aligned}
R\pi_*\cO_\pi(j)&=\Sym^jE^\vee &&(j\ge0),\\
R\pi_*\cO_\pi(j)&=0 &&(j\in[-h,-1]),\\
R\pi_*\cO_\pi(-\ell-d)&\simeq\big(\Sym^dE\otimes\det E\big)[-h] &&(d\ge0).
\end{aligned}
\end{equation}
For the third formula, write $\cO_\pi(-\ell-d)\cong G^\vee\otimes\omega_\pi$ with $G=\cO_\pi(d)\otimes\pi^*\det E^\vee$. Grothendieck duality gives $R\pi_*\cO_\pi(-\ell-d)\simeq(R\pi_*G)^\vee[-h]$, and $R\pi_*G=\Sym^dE^\vee\otimes\det E^\vee$.

\begin{proposition}\label{prop:gen-Gamma}
Let $1\le m\le h$. An admissible sequence $\mathbf p$ of maximal length normalized at $q=(1,-\ell)$ with $\Wing(\mathbf p)=\Gamma_m$ is full. The same holds if $\Wing(\mathbf p)=\rho(\Gamma_m)$.
\end{proposition}

\begin{proof}
Let $\cC=\gen{L(p):p\text{ a term of }\mathbf p}$. By Lemma~\ref{lem:basic-count}, the central columns $2,\dots,h$ are complete. The points $q=(1,-\ell)$ and $L_1,\dots,L_\ell$ form the complete column $1$. By Lemma~\ref{lem:complete}, $\cB_1,\dots,\cB_h\subseteq\cC$. By~\eqref{eq:rel-duality},
\begin{equation}\label{eq:Q-values}
\begin{gathered}
Q(L(0,0))=\cO_Y(-1)[-h],\qquad Q(L(R_t))=\cO_Y(t),\\
Q(L(D_z))=\big(\Sym^{z-1}E\otimes\cO_Y(z-1)\big)[-h].
\end{gathered}
\end{equation}
For $d\ge0$ put $H_d=(\Sym^dE)(d)$, so that $H_0=\cO_Y$. Let $\cD_0\subseteq\cD$ be the thick subcategory generated by the objects $Q(L(p))$, $p$ a term of $\mathbf p$. Since $\Wing(\mathbf p)=\Gamma_m$, it contains
\[
\cO_Y(-1),\ \cO_Y,\ \cO_Y(1),\ \cO_Y(m+1),\dots,\cO_Y(h),\qquad H_d\quad(1\le d<m).
\]
For $d\ge1$, the Euler sequence~\eqref{eq:euler} induces the exact sequence $0\to\Sym^dE\to\Sym^dU\otimes\cO_Y\to\Sym^{d-1}U\otimes\cO_Y(1)\to0$, and twisting by $\cO_Y(d)$ gives
\begin{equation}\label{eq:sym-euler}
0\fto H_d\fto\Sym^dU\otimes\cO_Y(d)\fto\Sym^{d-1}U\otimes\cO_Y(d+1)\fto0 .
\end{equation}
Start from $\cO_Y(1)\in\cD_0$ and use that $\cD_0$ is closed under direct summands. Then the sequences~\eqref{eq:sym-euler} for $d=1,\dots,m-1$ successively give $\cO_Y(2),\dots,\cO_Y(m)\in\cD_0$. Thus
\[
\cD_0\supseteq\gen{\cO_Y(-1),\cO_Y,\dots,\cO_Y(h)}=\cD
\]
by~\eqref{eq:SOD}, and Lemma~\ref{lem:Q} shows that $\mathbf p$ is full.

If $\Wing(\mathbf p)=\rho(\Gamma_m)$, then the sequence $\mathbf p^\tau$ of Lemma~\ref{lem:tau}(c) has wing $\Gamma_m$. It is therefore full, and so is $\mathbf p$.
\end{proof}

\begin{proposition}\label{prop:gen-T}
An admissible sequence $\mathbf p$ of maximal length normalized at $q=(1,-\ell)$ with $\Wing(\mathbf p)\in\{T_{00},T_{10},T_{01},T_{11}\}$ is full.
\end{proposition}

\begin{proof}
Let $\cC=\gen{L(p):p\text{ a term of }\mathbf p}$, and recall from Lemma~\ref{lem:basic-count} that the central columns are complete.

\emph{The wing $T_{00}$.} Rotate $\mathbf p$ once, replacing $0$ by $u=(\ell,\ell)$. Then $q,L_1,\dots,L_\ell$ form the complete column $1$, and $R_0,\dots,R_h,u$ form the complete column $\ell$. The rotated sequence is full by Lemma~\ref{lem:complete}, and so is $\mathbf p$ by Lemma~\ref{lem:operations}.

\emph{The wing $T_{01}$.} Column $1$ is complete, so $\cB_1,\dots,\cB_h\subseteq\cC$. By~\eqref{eq:Q-values} and~\eqref{eq:rel-duality}, the objects $Q(L(p))$ include, up to shift,
\[
\cO_Y(-1),\ \cO_Y,\ \cO_Y(1),\dots,\cO_Y(\ell-2),\qquad Q(L(F_h))=Q(L(\ell+1,-1))=E^\vee(-1)=T_Y(-2),
\]
where we used $E^\vee=T_Y(-1)$. Dualizing~\eqref{eq:euler} and twisting by $\cO_Y(-1)$ gives
\begin{equation}\label{eq:dual-euler}
0\fto\cO_Y(-2)\fto U^\vee\otimes\cO_Y(-1)\fto T_Y(-2)\fto0,
\end{equation}
so the thick subcategory generated by these images also contains $\cO_Y(-2)$. Hence it contains $\gen{\cO_Y(-2),\dots,\cO_Y(\ell-2)}=\cD$, and Lemma~\ref{lem:Q} applies.

\emph{The wing $T_{10}$.} The sequence $\mathbf p^\tau$ of Lemma~\ref{lem:tau}(c) has wing $\rho(T_{10})=T_{01}$. It is full by the previous case, and so is $\mathbf p$.

\emph{The wing $T_{11}$.} The complete central columns give $\cB_2,\dots,\cB_h\subseteq\cC$, but column $1$ now lacks $L_1$. The tautological surjection $\pi^*E^\vee\twoheadrightarrow\cO_\pi(1)$ has an exact Koszul complex. Tensoring it with $\cO_\pi(\ell)$ gives the exact complex
\begin{equation}\label{eq:koszul-X}
0\fto\pi^*\det E^\vee\fto\cO_\pi(1)\otimes\pi^*\Lambda^{h}E^\vee\fto\cdots\fto\cO_\pi(h)\otimes\pi^*E^\vee\fto\cO_\pi(\ell)\fto0,
\end{equation}
whose term of relative degree $i$ is $\cO_\pi(i)\otimes\pi^*\Lambda^{\ell-i}E^\vee\in\cB_i$. Its end terms are $\pi^*\det E^\vee=L(0,1)=L(D_1)$ and $\cO_\pi(\ell)=L(\ell,0)=L(R_0)$, both in $\cC$, and its terms of relative degree $2,\dots,h$ lie in $\cB_2,\dots,\cB_h\subseteq\cC$. Since $\Lambda^hE^\vee\simeq E\otimes\det E^\vee\simeq\Omega_Y(2)$, exactness gives
\[
\iota_1(\Omega_Y(2))=L(1,0)\otimes\pi^*\Omega_Y(2)\in\cC .
\]
The terms $q,L_2,\dots,L_\ell$ in column $1$ give $\iota_1(\cO_Y(-\ell))$ and $\iota_1(\cO_Y(j))$ for $2-\ell\le j\le0$. Splicing the exact sequences $0\to\Omega_Y^j(2)\to\Lambda^jU\otimes\cO_Y(2-j)\to\Omega_Y^{j-1}(2)\to0$ for $2\le j\le\ell$, which are twists of the exterior powers of~\eqref{eq:euler}, and using $\Omega_Y^\ell(2)=\cO_Y(1-\ell)\cong\Lambda^{\ell+1}U\otimes\cO_Y(1-\ell)$, we obtain the exact complex
\begin{equation}\label{eq:koszul-Y}
0\fto\Lambda^{\ell+1}U\otimes\cO_Y(1-\ell)\fto\Lambda^{\ell}U\otimes\cO_Y(2-\ell)\fto\cdots\fto\Lambda^2U\otimes\cO_Y\fto\Omega_Y(2)\fto0 .
\end{equation}
Applying $\iota_1$, every term except the first lies in $\cC$; hence $\iota_1(\cO_Y(1-\ell))=L(1,1-\ell)\in\cC$. Column $1$ now contains $\iota_1(\cO_Y(j))$ for all $-\ell\le j\le0$, so $\cB_1\subseteq\cC$. The terms $0,R_0,\dots,R_{h-1},F_h$ are the same as in the case $T_{01}$, so~\eqref{eq:dual-euler} and Lemma~\ref{lem:Q} finish the proof.
\end{proof}

\section{Proof of the descent criterion}\label{sec:proof}

\begin{proof}[Proof of Theorem~\ref{thm:descent}]
Let $\mathbf p$ be an admissible sequence of maximal length with a non-orthogonal descent. By Lemma~\ref{lem:operations}, we may assume that $\mathbf p$ is normalized at $q$ as in~\eqref{eq:normalized}. By Lemma~\ref{lem:descent-shapes}, either $q$ has shape~\ref{shape:diag}, or $q=(s,-k)$ with $1\le s\le h$, $k\ge\ell$ and $k\ne\ell+s$. Shape~\ref{shape:diag} is impossible by Proposition~\ref{prop:diagonal}, and $s\ge2$ is impossible by Corollary~\ref{cor:no-strip}. Hence $q=(1,-k)$ with $k=\ell$ or $k\ge\ell+2$. If $k\ge\ell+2$, then $\mathbf p$ is full by Proposition~\ref{prop:basic-large}. If $k=\ell$, then $\Wing(\mathbf p)$ is an admissibly orderable subset of $W_q$ of size $2\ell$ by Lemma~\ref{lem:basic-count}, hence one of the sets listed in Theorem~\ref{thm:classification}. In each case $\mathbf p$ is full by Propositions~\ref{prop:gen-Gamma} and~\ref{prop:gen-T}.
\end{proof}

This completes the proof of Theorem~\ref{thm:main}, and Corollary~\ref{cor:main} follows by combining it with~\cite[Theorem~C]{AHW25}.

\begin{remark}\label{rem:l2}
The hypothesis $\ell\ge3$ is used at several places, for instance in the incompatibility of $F_1-L_1=(2\ell-2,0)$ in Lemma~\ref{lem:basic-count} and in the case analysis of Proposition~\ref{prop:s2}. For $\ell=2$ the conclusion of Theorem~\ref{thm:main} holds by~\cite[Theorem~C]{AHW25}, which is proved there by different methods; we have not attempted to adapt the case analysis above to $\ell=2$.
\end{remark}

\end{document}